\documentclass[10pt,reqno]{amsart}

\usepackage{amssymb}
\usepackage{amsfonts}
\usepackage{amsmath}
\usepackage{latexsym}
\usepackage[dvips]{graphicx}
\usepackage{psfrag}
\usepackage{amsopn}
\usepackage{amstext}
\usepackage{amscd}
\usepackage{amsthm}
\usepackage[centertags,leqno]{amsmath}

\newtheorem{theorem}{Theorem}
\newtheorem{proposition}[theorem]{Proposition}
\newtheorem{corollary}[theorem]{Corollary}
\newtheorem{lemma}[theorem]{Lemma}
\theoremstyle{definition}
\newtheorem{definition}[theorem]{Definition}
\theoremstyle{remark}
\newtheorem{remark}[theorem]{Remark}

\title{SELF--SCALED BARRIERS FOR IRREDUCIBLE SYMMETRIC CONES}
\author{Raphael A.\ Hauser, Yongdo Lim}

\DeclareMathOperator{\e}{e} 
\DeclareMathOperator{\Int}{int} \DeclareMathOperator{\GL}{GL}
\DeclareMathOperator{\tr}{tr} \DeclareMathOperator{\Aut}{Aut}
 \DeclareMathOperator{\Orth}{O} 
\DeclareMathOperator{\SO}{SO}
\DeclareMathOperator{\id}{id} \DeclareMathOperator{\cst}{cst}
\DeclareMathOperator{\Iso}{Iso} \DeclareMathOperator{\T}{T}
\DeclareMathOperator{\I}{I}

\newcommand{\Half}{\ensuremath{{1/2}}}
\newcommand{\comment}[1]{}
\newcommand{\Ck}[1]{\ensuremath{\mathcal{C}^{#1}}}

\newcommand{\RN}{\ensuremath{\mathbb{R}}}

\newcommand{\Sym}{\ensuremath{{\mathrm Sym}(n,\RN)}}
\newcommand{\Spd}{\ensuremath{\Sym^{+}}}

\begin{document}

\thispagestyle{empty}

\vspace*{45pt}

\begin{center}
{\Large UNIVERSITY OF CAMBRIDGE}

\vspace*{6mm}
{\large Numerical Analysis Reports}

\vspace*{30mm}
\parbox[t]{105mm}{\large\bf\begin{center}
SELF--SCALED BARRIERS FOR IRREDUCIBLE SYMMETRIC CONES 
\end{center}}

\vspace*{15mm}
{\bf Raphael Hauser and Yongdo Lim}

\vspace*{45mm}
DAMTP 2001/NA04\\ 
April 2001

\vspace*{-330pt}
\begin{picture}(370,340)(0,0)
\put (3,3) {\line(1,0){364}}
\put (3,3) {\line(0,1){334}}
\put (367,3) {\line(0,1){334}}
\put (3,337) {\line(1,0){364}}
\thicklines
\put (0,0) {\line(1,0){370}}
\put (0,0) {\line(0,1){340}}
\put (370,0) {\line(0,1){340}}
\put (0,340) {\line(1,0){370}}
\end{picture}

\vfill
{\large Department of Applied Mathematics and Theoretical
Physics\\Silver Street\\Cambridge\\England CB3 9EW}
\end{center}

\thispagestyle{empty}

\thispagestyle{empty}
\newpage

\hspace{3cm}

\vfill
\thispagestyle{empty}

\newpage
\setcounter{page}{1}

\maketitle

\vspace{.7cm}

\begin{center}
April 2, 2001
\end{center}

\vspace{1cm}

\begin{abstract}
\noindent Self--scaled barrier functions are fundamental objects
in the theory of interior--point methods for linear optimization
over symmetric cones, of which linear and semidefinite programming
are special cases. We are classifying all self--scaled barriers
over irreducible symmetric cones and show that these functions are
merely homothetic transformations of the universal barrier
function. Together with a decomposition theorem for self--scaled
barriers this concludes the algebraic classification theory of
these functions. After introducing the reader to the concepts
relevant to the problem and tracing the history of the subject, we
start by deriving our result from first principles in the
important special case of semidefinite programming. We then
generalise these arguments to irreducible symmetric cones by
invoking results from the theory of Euclidean Jordan algebras.
\end{abstract}

\vspace{2cm}

{\tiny
\noindent
{\bf Key Words}\\
Semidefinite programming, self--scaled barrier functions,
interior--point methods, symmetric cones, Euclidean Jordan
algebras.\\
{\bf AMS 1991 Subject Classification}\\
Primary 90C25, 52A41, 90C60. Secondary 90C05, 90C20.\\
{\bf Contact Information and Credits}\\
Raphael Hauser, Department of Applied Mathematics and Theoretical
Physics, University of Cambridge, Silver Street, Cambridge CB3
9EW, England. {\em rah48@damtp.cam.ac.uk}. Research supported in
part by the Norwegian Research Council through project No.\
127582/410 ``Synode II'', by the Engineering and Physical Sciences
Research Council of the UK under grant No.\ GR/M30975,
and by NSERC of Canada grants of J.\ Borwein and P.\ Borwein.\\
Yongdo Lim, Department of Mathematics, Kyungpook National
University, Taegu 702--701, Korea. {\em ylim@knu.ac.kr}. 
Research supported in part by  the Basic Research Program of the 
Korea Science and  Engineering Foundation through project  
No. 2000-1-10100-007-3.}

\newpage

\section{Introduction}
{\em Self--scaled barriers} are a special class of
self--concordant barrier functions \cite{int:Nesterov5} introduced
by Nesterov and Todd \cite{int:Nesterov10} for the purpose of
extending long--step primal--dual symmetric interior--point
methods from linear and semidefinite programming to more general
convex optimization problems
\cite{int:Nesterov10,int:Nesterov12}, see Definition~\ref{selfscaled} 
below.

The domain of definition of a self--scaled barrier $F$ is an proper 
open convex cone $\Omega$ lying in a real Euclidean space 
$\bigl(V,\langle\cdot,\cdot\rangle\bigr)$. By abuse of language 
one often refers to $F$ as a self--scaled barrier for the 
topological closure $\bar{\Omega}$ of $\Omega$. Not every proper 
open convex cone $\Omega$ allows a self--scaled barrier, and for 
those who do $\bar{\Omega}$ is called a {\em self--scaled cone} in 
the terminology of Nesterov and Todd 
\cite{int:Nesterov10}. G\"uler \cite{int:Gueler10} found that the
family of interiors of self--scaled cones is identical to the set of
{\em symmetric cones} studied in the theory of Euclidean Jordan
algebras. Due to this discovery, Jordan algebra theory became an 
important analytic tool in the theory of semidefinite programming 
and its natural generalisation, self--scaled programming.

Self--scaled barriers were studied by Nesterov--Todd
\cite{int:Nesterov10,int:Nesterov12}, G\"uler \cite{int:Gueler10},
G\"uler--Tun{\c c}el \cite{Levent} (see p.\ 124 and related
material), Hauser \cite{Hauserthesis} and others. Though an
axiomatic theory of these functions exists, the only known
examples are trivially related to the characteristic function of
the cone $\Omega$,
\begin{equation}\label{eq:charfunction}
\varphi_{\Omega}(x):=\int_{\Omega^\sharp}\e^{-\langle x,s\rangle}ds,
\end{equation}
where $\Omega^{\sharp}:=\{s\in V:\langle x,s\rangle>0\;\forall
x\in\Omega\}$ is the polar of $\Omega$ with respect to
$\langle\cdot,\cdot\rangle$. This was first discovered by G\"uler 
who showed that the universal barrier $U$ \cite{int:Nesterov5} of a 
symmetric cone is self--scaled and is a homothetic transformation
$U=c_1\ln\varphi_{\Omega}+c_2$ of the characteristic function
$\varphi_{\Omega}$, where $c_1\geq 1$ and $c_2$ are constants. More
generally, every symmetric cone $\Omega$ has a decomposition, 
unique up to indexing, into a direct sum of irreducible symmetric 
cones
\[\Omega=\Omega_1\oplus\dots\oplus\Omega_m,\]
where the $\Omega_i$ lie in subspaces $V_i\subseteq V$ decomposing $V$
into a direct sum $V=\oplus_{i=1}^m V_i$. The {\em irreducible summands} 
$\Omega_i$ can be classified into five different types, see \cite{FK}
and references therein. All known self--scaled barrier functions for 
$\Omega$ are of the form
\begin{equation}\label{isoform}
H=c_0+\bigoplus_{i=1}^m c_i\ln\varphi_{\Omega_i},
\end{equation}
where $c_1,\dots,c_m\geq 1$. We use the direct sum notation for $H$ in
Equation \eqref{isoform} and elsewhere to indicate that each
$x=\oplus_{i=1}^m x_i\in\oplus_{i=1}^m\Omega_i=\Omega$ is
mapped to $H(x)=c_0+\sum_{i=1}^m c_i\ln\varphi_{\Omega_i}(x_i)$. 
It is also well--known that each function of the form \eqref{isoform} 
is a self--scaled barrier for $\Omega$.

The natural question arises as to whether {\em all} self--scaled
barrier functions are of the form \eqref{isoform}. Early dents
into this question were made by G\"uler and Tun{\c c}el when
considering {\em invariant barriers}, see \cite{Levent} page 124
and related material. In a chapter of his thesis
\cite{Hauserthesis} and in a subsequent report \cite{Hauser99},
Hauser showed that any self--scaled barrier $H$ over a symmetric cone 
$\Omega$ decomposes into a direct sum $H=\oplus_{i=1}^m H_i$ of
self--scaled barriers $H_i$ over the irreducible components
$\Omega_i$, and that any {\em isotropic}, i.e., rotationally invariant,
self--scaled barrier $H_i$ on $\Omega_i$ is of the form
$c_1\ln\varphi_{\Omega_i}+c_2$ with $c_1\geq 1$. Hauser also observed 
that any self--scaled barrier $H$ on $\Omega$ is invariant under a 
rich class of rotations of $\Omega$, i.e., elements of
$\Orth(\Omega)$, see \eqref{eq:orthogonalgroup} below, where these 
particular rotations are defined in terms of the Hessians of $H$, see 
Lemma 2.2.19 \cite{Hauserthesis}. Suspecting that in the case where 
$\Omega$ is irreducible this family of rotations is rich enough to generate
all of $\Orth(\Omega)$, Hauser \cite{Hauserthesis,Hauser99} conjectured
that all self--scaled barriers on irreducible symmetric cones are
isotropic. This conjecture, whose correctness to prove is the
primary objective of this paper, shows that all self--scaled
barriers are indeed of the form
\eqref{isoform} and concludes their algebraic classification.

In a second report \cite{Hauser00} Hauser showed the correctness
of the isotropy conjecture for the special case of the cone of
positive semidefinite symmetric matrices. The proof follows
exactly the path suggested by Lemma 2.2.19 \cite{Hauserthesis}, and as
outlined above. Shortly after Hauser's report \cite{Hauser00} was
announced, Lim \cite{LimPC} generalised Hauser's arguments to
general irreducible symmetric cones and settled the isotropy
conjecture. This report is a joint publication consisting of a 
revision of Hauser's report \cite{Hauser00} and of Lim's 
generalisation \cite{LimPC}.

Subsequently, both Schmieta~\cite{Schmieta} and
G\"uler~\cite{Guler00} independently of each other and
independently of Lim also proved the isotropy conjecture.
Schmieta's report \cite{Schmieta} was the first publication where
the full classification result became available. 
G\"uler's approach \cite{Guler00} was later incorporated in a 
joint publication with Hauser \cite{HauserGuler} which started as a 
major revision of the report \cite{Hauser99}.

It is interesting to note that, though the approaches of 
Lim, G\"uler and Schmieta differ in important details, all
three involve two key mechanisms: The so--called 
{\em fundamental formula} on the one hand, see \eqref{eq:fundamental} 
below, and Koecher's Theorem 4.9 (b) \cite{Koe} on the other hand. 
Already Hauser's approach \cite{Hauser00} to solving the special case 
of the positive semidefinite cone was based on fundamentally the same 
ideas, as his Proposition 3.3 was essentially an independent 
rediscovery of Koecher's theorem in this particular case (c.f.\ 
Corollary 4.3 \cite{Hauser00}).

The main part of this paper is organised in two sections
addressing different communities: Section~\ref{SDP} treats the
case of the positive semidefinite cone only. Readers interested in
semidefinite programming and lacking a background in Jordan
algebra theory will find it easy to read this section, the results
being derived from first principles. All the qualitative features
of the general approach already appear in the restricted
framework, and it is possible to understand some of the essential
ideas behind Koecher's theorem by reading the proof of
Proposition~\ref{generate}. Section~\ref{ISC} on the other hand treats the
general case and addresses primarily readers with a background in
Jordan algebra theory. It is possible to understand this section
without prior lecture of Section~\ref{SDP}.

We conclude this section by giving the essential definitions and
identities that form the basis of the theory of self--scaled
barriers. Recall that we introduced the notation $\Omega$ and $V$ 
above. The set of vector space automorphisms of $V$ that leave 
$\Omega$ invariant is 
called the {\em symmetry group} of $\Omega$, and we denote it by
\begin{equation*}
G(\Omega):=\bigl\{\theta\in\GL(V):\theta\Omega=\Omega\bigr\}.
\end{equation*}
The inner product on $V$ defines a notion of adjoint of an 
endomorphism and hence an orthogonal group 
$\Orth(V)=\bigl\{\theta\in \GL(V): \theta^*=\theta^{-1}\bigr\}$. 
The subgroup 
\begin{equation}\label{eq:orthogonalgroup}
\Orth(\Omega):=G(\Omega)\cap\Orth(V)
\end{equation}
is called the {\em orthogonal group} of $\Omega$. 

\begin{definition}[\cite{int:Nesterov10}]\label{selfscaled}
A $\nu$\nobreakdash--self--concordant logarithmically homogeneous
\cite{int:Nesterov5} barrier functional $H\in\Ck{3}(\Omega,\RN)$
is said to be {\em self--scaled} if the following conditions are
satisfied:
\begin{itemize}
\item[(a)] $H''(w)x\in\Omega^{\sharp}$ for all
$x,w\in\Omega$ and
\item[(b)] $H_\sharp\bigl(H''(w)x\bigr)=H(x)-2H(w)-\nu$
for all $x,w\in\Omega$.
\end{itemize}
\end{definition}

The function $H_\sharp:s\mapsto\max\bigl\{-\langle
x,s\rangle-H(x):x\in\Omega\bigr\}$ is a self--scaled barrier
defined on $\Omega^{\sharp}$ \cite{int:Nesterov10}. It is assumed 
in the definition of a self--concordant function, see 
\cite{int:Nesterov5}, that the Hessians $H''(x)$ are non--singular 
for all $x\in\Omega$. The next theorem is a compilation of several 
separate results of Nesterov and Todd \cite{int:Nesterov10}:

\begin{theorem}[\cite{int:Nesterov10}]\label{thm2}
Let $H\in\Ck{3}(\Omega,\RN)$ be self--scaled and $x\in\Omega,
s\in\Omega^\sharp$. Then there exists a unique \emph{scaling
point} $w_H(x,s)\in\Omega$ such that $s=H''\bigl(w(x,s)\bigr)x$. 
Furthermore, the following properties hold:
\begin{itemize}
\item[(a)] $H''(w)\in\Iso\bigl(\Omega,\Omega^\sharp\bigr)
\quad\forall w\in\Omega,$
\item[(b)] $H''(x)=H''\bigl(w(x,s)\bigr)\circ H_{\sharp}''(s)\circ 
H''\bigl(w(x,s)\bigr).$\qed
\end{itemize}
\end{theorem}

It is customary to change the inner product
$\langle\cdot,\cdot\rangle$ to
$\langle\cdot,\cdot\rangle_f:=\langle H''(f)\cdot,\cdot\rangle$
where $f$ is a fixed element in $\Omega$. We will always assume
that $\langle\cdot,\cdot\rangle$ is already of this kind, i.e.,
that there exists an element $f\in\Omega$ such that $H''(f)=\id_{\Omega}$
is the identity when the Hessian is computed with respect to this
inner product. Under this assumption it is easy to show that 
$\Omega^{\sharp}=\Omega$ and $H_\sharp=H+\cst$. 
Hence, in this framework we can reformulate Parts $(a)$ and $(b)$ 
of Theorem~\ref{thm2} as follows:
\begin{align}
H''(w)&\in G(\Omega),\label{eq:thm2i}\\
H''(x)&=H''\bigl(w_H(x,s)\bigr)\circ H''(s)\circ
H''\bigl(w_H(x,s)\bigr)\label{eq:thm2ii}.
\end{align}

Theorem~\ref{thm2} reveals that $\Omega$ can allow a self--scaled
barrier only if it is a {\em symmetric cone}, c.f.\
Section~\ref{ISC}. As mentioned above, this fact was first 
observed by G\"uler \cite{int:Gueler10} who showed that the
relation also holds in the opposite direction.

Equation \eqref{eq:thm2ii} is a reformulation of an identity which
is called {\em fundamental formula} in Jordan algebra theory, see
Equation~\eqref{eq:fundamental} below. This identity is one of the keys to
proving the isotropy conjecture, as it allows to express all the
Hessians of $H$ in terms of the Hessian $H''(e)$ at a single
element $e\in\Omega$ and of the Hessians of the standard
logarithmic barrier function
\begin{equation}\label{standard}
F(x)=\ln\varphi_{\Omega}(x),
\end{equation}
see Equation~\eqref{eq:charfunction},  which is self--scaled 
\cite{int:Gueler10}. Rothaus \cite{Rothaus} proved the following
result that will be important for our purposes:

\begin{theorem}[\cite{Rothaus}]\label{polar}
For every $\theta\in G(\Omega)$ there exists a unique
$\omega\in\Orth(\Omega)$ and a unique $w\in\Int(\Omega)$ such that $\theta$
has the {\em polar decomposition} $\theta=\omega\circ F''(w)$. \qed
\end{theorem}

Since $\Orth(\Omega)\subset\Orth(V)$ and $F''(w)$ is a self--adjoint
positive definite automorphism of $V$, Rothaus' polar decomposition is
identical to Cartan's polar decomposition. Theorem~\ref{polar} shows 
the non--trivial fact that both factors lie in $G(\Omega)$. The
uniqueness of Cartan's polar decomposition trivially implies the
following lemma which will be useful in later sections:

\begin{lemma}\label{self-adjoint}
The set of self--adjoint positive definite automorphisms of $V$
that preserve $\Omega$ coincides with
$\bigl\{F''(w):w\in\Omega\bigr\}$. \qed
\end{lemma}

\section{Self--Scaled Barriers for Semidefinite Programming}\label{SDP}

This section is limited to semidefinite programming and provides 
readers who are unfamiliar with Jordan algebra terminology access
to the main ideas behind the mechanism that forces self--scaled 
barriers to be essentially unique. 

In this framework it is customary to write variable names with
capitalised letters. $V$ is the space $\Sym$ of symmetric $n\times
n$ matrices endowed with the trace inner product $\langle
X,S\rangle=\tr\bigl(X^{*}S\bigr)$ which corresponds to the
Frobenius norm. $\Omega$ is the cone $\Spd$ of positive definite 
symmetric $n\times n$ matrices. The following identities hold for 
the standard logarithmic barrier function $F(X)=-\ln\det(X)$:
\begin{align}
\label{eq:SDPHessian}&F^{\prime\prime}(X)A
=X^{-1}AX^{-1}\qquad\forall A\in\Sym\quad\text{and}\\
\label{eq:SDPScaling}&W_F(X,S)=S^{-\Half}\bigl(S^{\Half}X
S^{\Half}\bigr)^{\Half}S^{-\Half}.
\end{align}

Let us assume that $H$ is an arbitrary self--scaled barrier
function for $\Spd$. Applying Equation \eqref{eq:thm2ii} to both
$H$ and $F$, we can derive a series of expressions that will allow
us to relate Hessians of $H$ to Hessians of $F$. Since $H''(X)$ is
a self--adjoint positive definite automorphism of $\Spd$ it follows
from Lemma \ref{self-adjoint} that there exists a well--defined
mapping $\Upsilon:\Spd\rightarrow\Spd$ such that 
\begin{equation}\label{eq:Yscaling}
H''(X)=F''\bigl(\Upsilon(X)\bigr)
\end{equation}
for all $X\in\Spd$. We claim that $\Upsilon$ is a scalar function,
i.e., there exists a number $\lambda>0$ such that 
$\Upsilon=\lambda\cdot\id_{\Sym}$. The proof of this claim is going to
occupy us until Corollary~\ref{hessians}.

Let $W_F=W_F(X,S)$ and $W_H=W_H(X,S)$ denote the scaling points of 
$X,S\in\Spd$ defined by $F$ and $H$, see Theorem~\ref{thm2}. 
Equation \eqref{eq:Yscaling} implies that 
$F''\bigl(\Upsilon(W_H)\bigr)X=S=F''\bigl(W_F\bigr)X$, and it
follows from the uniqueness part of Theorem \ref{thm2} that
$\Upsilon\bigl(W_H\bigr)=W_F$. Therefore, Equation
\eqref{eq:thm2ii} applied to $H$ shows that
\begin{align}
H''(X)&=F''(W_F)\circ H''(S)\circ F''(W_F),\quad\text{and}
\label{eq:FHscaling}\\
\Upsilon(X)&=W_F\Upsilon(S)W_F.\label{eq:ypsscale}
\end{align}
Note that $F''(\I)=\id_{\Sym}$. Therefore, Equation \eqref{eq:thm2ii}
applied to $F$ and $S=\I$ shows that
$F''\bigl(W_F(X,\I)\bigr)=F''(X)^{-\Half}$. Using this fact in 
conjunction with Equations \eqref{eq:SDPHessian} and 
\eqref{eq:FHscaling} we get 
\begin{equation}\label{eq:XIYps}
\Upsilon(X)=X^{\Half}\Upsilon(\I)X^{\Half}\quad\forall X\in\Spd.
\end{equation}
Therefore,
\begin{align*}
X^{\Half}\Upsilon(I)X^{\Half}
\stackrel{\eqref{eq:XIYps},\eqref{eq:ypsscale}}{=}&
W_F(X,S)\Upsilon(S)W_F(X,S)\\
\stackrel{\eqref{eq:XIYps}}{=}&
W_F(X,S)S^{\Half}\Upsilon(I)S^{\Half}W_F(X,S),
\end{align*}
and by virtue of \eqref{eq:SDPScaling} this implies that
\begin{equation}\label{eq:star}
\Upsilon(I)=X^{-\Half}S^{-\Half}\bigl(S^{\Half}XS^{\Half}\bigr)^{\Half}
\Upsilon(I)S^{-\Half}\bigl(S^{\Half}XS^{\Half}\bigr)^{\Half}X^{-\Half}
\end{equation}
for all $X,S\in\Spd$. Clearly, this condition is equivalent to
\begin{equation}\label{eq:fdmt}
\Upsilon(I)=N^{-1}\bigl(NN^{\T}\bigr)^{\Half}
\Upsilon(I)\bigl(NN^{\T}\bigr)^{\Half}N^{-1}\quad\forall
N\in{\mathcal K},
\end{equation}
where ${\mathcal K}:=\bigl\{XS:X,S\in\Spd\bigr\}$ is the set of
$n\times n$ matrices that can be written as the product of two
symmetric positive definite matrices. The following characterisation 
of this set is due to Mike Todd \cite{Mike Todd}:

\begin{lemma}\label{nondefect}
${\mathcal K}$ coincides with the set of non--defective $n\times n$
matrices with real coefficients, all of whose eigenvalues are
strictly positive real numbers.
\end{lemma}

\begin{proof}
If $N=XS$, then 
\begin{equation*}
N=X^{\Half}\bigl(X^{\Half}SX^{\Half}\bigr)X^{-\Half}
=X^{\Half}\bigl(QDQ^{\T}\bigr)X^{-\Half}=PDP^{-1},
\end{equation*}
where $QDQ^{\T}$ is the spectral decomposition of the symmetric
positive definite matrix $X^{\Half}SX^{\Half}$, and
$P:=X^{\Half}Q$. This gives the eigenvalue decomposition of $N$,
with eigenvalues the positive entries of $D$ and eigenvectors the
columns of $P$. Conversely, suppose $N=PDP^{-1}$, where $D$ is
diagonal with positive diagonal entries. Then we can write
$N=\bigl(PP^{\T}\bigr)\bigl(P^{-T}DP^{-1}\bigr)=:XS$.
\end{proof}

Note that $N^{-1}(NN^{\T})^{\Half}\in\SO(n)$ for all
$N\in{\mathcal K}$. In the next proposition we characterise the
closed subgroup of $\SO(n)$ generated by matrices of this form.
This constitutes the mathematical core mechanism of our proof and 
is a close relative of Koecher's Theorem 4.9 (b) \cite{Koe}.

\begin{proposition}\label{generate}
The closed subgroup of $\SO(n)$ generated by the set of orthogonal
matrices of the form $N^{-1}(NN^{\T})^{\Half}$ with 
$N\in{\mathcal K}$ coincides with $\SO(n)$.
\end{proposition}

\begin{proof}
Let this closed subgroup be denoted by ${\mathcal G}$, and let 
${\mathfrak g}$ be its Lie Algebra. Since $\SO(n)$ is connected, it suffices
to show that ${\mathfrak g}=\mathfrak{so}(n)$, or in other words
that the tangent space of the manifold ${\mathcal G}$ at the point 
$\I\in{\mathcal G}$ coincides with the set of $n\times n$ skew--symmetric
matrices over the reals. In fact, use of the exponential mapping
$\exp:\mathfrak{so}(n)\rightarrow\SO(n)$ 
shows that ${\mathcal G}$ and $\SO(n)$ share a neighbourhood 
${\mathcal V}$ of $\I$, and parallel transport by left trivialisation 
shows that ${\mathcal G}$ 
and $\SO(n)$ share the neighbourhood $g{\mathcal V}$ of any element
$g\in{\mathcal G}$. Therefore, ${\mathcal G}$ is both open and closed
as a subset of $\SO(n)$, and since $\SO(n)$ is connected, the result 
follows.

It remains to show that ${\mathfrak g}=\mathfrak{so}(n)$:
Let $\Delta\in{\mathcal K}$ have eigenvalues
$\lambda_1,\dots\lambda_n>0$. Then $\Delta(t):=I+t\Delta
\in{\mathcal K}$ for all $t>0$, since the $n$ linearly independent
eigenvectors of $\Delta$ are also eigenvectors of $\Delta(t)$ and
correspond to the strictly positive eigenvalues
$t\lambda_i(t)+1>0$, $(i=1,\dots,n)$. The Neumann--series
development of $\Delta(t)^{-1}$ shows that
\begin{equation}\label{neumann}
\Delta(t)^{-1}=I-t\Delta+O(t^2).
\end{equation}
Upon taking squares on both sides of the ansatz
\begin{equation}\label{ansatz}
(\Delta(t)\Delta(t)^{\T})^{\Half}=I+tD+O(t^2)
\end{equation}
we get $I+t(\Delta+\Delta^{\T})+O(t^2)=I+2tD+O(t^2)$, and hence
$D=\Half(\Delta+\Delta^{\T})$. Equations \eqref{neumann} and
\eqref{ansatz} thus yield the identity
\begin{align*}
\Delta(t)^{-1}\bigl(\Delta(t)\Delta(t)^{\T}\bigr)^{\Half}&=
\bigl(I-t\Delta\bigr)\bigl(I+\frac{t}{2}(\Delta+\Delta^{\T})\bigr)+O(t^2)\\
&=I+\frac{t}{2}(\Delta^{\T}-\Delta)+O(t^2),
\end{align*}
and this shows that $\Delta^{\T}-\Delta\in{\mathfrak g}$ for all
$\Delta\in{\mathcal K}$. Clearly, we have
$\Delta^{\T}-\Delta\in\mathfrak{so}(n)$ as expected. Now, for
$P:=\bigl(\begin{smallmatrix}1&0\\2&1\end{smallmatrix}\bigr)$ and
$D:=\bigl(\begin{smallmatrix}.5&0\\0&1\end{smallmatrix}\bigr)$ we
get $P^{-1}DP=
\bigl(\begin{smallmatrix}.5&1\\0&1\end{smallmatrix}\bigr)$. Hence,
$(P^{-1}DP)^{\T}-P^{-1}DP=
\bigl(\begin{smallmatrix}0&1\\-1&0\end{smallmatrix}\bigr)$. Let
$P_{ij}$ be the permutation matrix that permutes the $i^{th}$ and
$j^{th}$ variables, and let $E_{ij}^{-}:=e_ie_j^{\T}-e_je_i^{\T}$,
where $e_i$ is the $i$--th element in the canonical basis of
$\RN^n$. Then consider
\begin{equation*}
\Delta:=\bigl(P_{1i}P_{2j}
\bigl(\begin{smallmatrix}P^{-1}&\\&\I\end{smallmatrix}\bigr)\bigr)
\bigl(\begin{smallmatrix}D&\\&\I\end{smallmatrix}\bigr)
\bigl(\bigl(\begin{smallmatrix}P&\\&\I\end{smallmatrix}\bigr)
P_{j2}P_{i1}\bigr).
\end{equation*}
Clearly, $\Delta\in{\mathcal K}$ and
$\Delta^{\T}-\Delta=E_{ij}^{-}$. But since
$\bigl\{E_{ij}^{-}:i,j\in\{1,\dots,n\}\bigr\}$ forms a basis of
$\mathfrak{so}$ we find that the elements of ${\mathfrak g}$ span
this whole space. This shows the claim.
\end{proof}

\begin{corollary}\label{hessians}
There exists a positive constant $\lambda>0$ such that
$H^{\prime\prime}(X)=\lambda F^{\prime\prime}(X)$ for all
$X\in\Spd$.
\end{corollary}

\begin{proof}
The invariance property \eqref{eq:fdmt} is clearly preserved when
taking compositions and limits. Hence, Lemma \ref{nondefect}
implies that the symmetric positive definite matrix $\Upsilon(\I)$
satisfies the condition $\Upsilon(I)=\Omega\Upsilon(\I)\Omega^{\T}$ 
for all $\Omega\in\SO(n)$. It is a trivial matter to prove that this
forces $\Upsilon(\I)$ to be a scalar, and the result now follows
from \eqref{eq:XIYps}.
\end{proof}

\begin{theorem}\label{thm:mainSDP}
If $H$ is a self--scaled barrier functional for the cone of
symmetric positive semidefinite $n\times n$ matrices then there
exist constants $c_1>0$ and $c_0\in\RN$ such that
\begin{equation*}
H:X\mapsto c_0-c_1\ln\det(X)\quad\forall X\in\Spd.
\end{equation*}
\end{theorem}

\begin{proof}
It follows from Corollary~\ref{hessians} and the fundamental theorem
of differential and integral calculus that $H$ is of the form
$c_1F+\varphi+c_0$, where $c_1=\lambda>0$, $c_0\in\RN$ and
$\varphi\in\Sym^\sharp$ is a linear form on $\Sym$, i.e.\ there
exists a symmetric matrix $Y\in\Sym$ such that
$\varphi:X\mapsto\tr(Y^{\T}X)$ for all $X\in\Sym$. One of the
conditions in the definition of a $\nu$--self--concordant barrier
functional $B$ for a convex open domain $D$ is that the length of
the Newton step $B''(x)^{-1}[-B'(x)]$ at $x\in D$, measured in the
Riemannian metric $\|.\|_x$ defined by $B''(x)$ be uniformly
bounded by $\nu^{\Half}$, see e.g.\
\cite{int:Nesterov5,int:Nesterov10,Renegar:book}, i.e.,
\begin{align*}
\|B''(x)^{-1}[-B'(x)]\|_x^2&:= \bigl\langle
B''(x)\bigl[-B''(x)^{-1}[B'(x)]\bigr],
-B''(x)^{-1}B'(x)\bigr\rangle\\
&=\bigl\langle B'(x),\bigl(B''(x)\bigr)^{-1}[B'(x)]\bigr\rangle
\leq \nu.
\end{align*}
In particular, in the case of $H$ this means that
\begin{align*}
\nu&\geq\|H'(X)\|_X^2=\|Y-\lambda X^{-1}\|_X^2
\stackrel{C.S.}{\geq}
\bigl(\|Y\|_X-\|\lambda X^{-1}\|_X\bigr)^2\\
&\stackrel{{\mathrm Cor.}\ref{hessians}}{=} \Bigl(\tr(\lambda^{-1}XYX\cdot
Y^{\T})^{\Half}- \bigl(\lambda^{-1}X(\lambda X^{-1})X\cdot(\lambda
X^{-1})^{\T}\bigr)^{\Half}
\Bigr)^2\\
&=\lambda^{-1}\Bigl(\bigl(\tr\bigl[(XY)^2\bigr]\bigr)^{\Half}-
\lambda n^{\Half}\Bigr)^2
\end{align*}
for all $X\in\Spd$. But clearly, this implies that $Y=0$.
\end{proof}

\begin{remark} It is possible to prove Theorem~\ref{thm:mainSDP} by
invoking Lemma 2.2.19 \cite{Hauserthesis}, see \cite{Hauser00} for
details. This approach is interesting, since it follows the path
traced by the intuition that first led to the isotropy conjecture.
However, the proof we gave above fits better into the mainstream 
literature and is therefore easier to understand.
\end{remark}

\section{Self--Scaled Barriers for Irreducible Symmetric Cones}\label{ISC}

An open convex cone $\Omega$ in a real Euclidean 
space $V$ that is self--dual with respect to the given inner product 
and is homogeneous in the sense that the group 
$$G(\Omega):=\{g\in {\mathrm{GL}}(V): g(\Omega)=\Omega\}$$ acts
transitively on $\Omega$ is called a \textit{symmetric cone}. 
The theory of symmetric cones is 
closely tied to that of Euclidean Jordan algebras. We recall 
certain basic notions and well--known facts concerning Jordan 
algebras from the book \cite{FK} by J. Faraut and A. Kor\'anyi.

A {\em Jordan algebra} $V$ over  the  field $\Bbb R$ or $\Bbb C$ 
is a commutative algebra satisfying $x^{2}(xy)=x(x^{2}y)$ for all
$x,y\in V.$ We also assume the existence of a multiplicative
identity $e.$ Denote by $L$ the left translation  $L(x)y=xy,$ and
$P$ by the quadratic representation $P(x)=2L(x)^{2}-L(x^{2})$ for
$x\in V.$ An alternate statement of the Jordan algebra law is
$(xy)x^{2}=x(yx^{2}),$ a weak associativity condition that is
strong enough to ensure that the subalgebra generated by $\{e,x\}$
in $V$ is associative. An element $x\in V$ is said to be
\textit{invertible} if there exists an element $y$ in the
subalgebra generated by $x$ and $e$ such that $xy=e.$ It is known
that an element $x$ in $V$  is invertible if and only if $P(x)$ is
invertible. In this case, $P(x)^{-1}=P(x^{-1}).$ If $x$ and $y$ are 
invertible, then $P(x)y$ is invertible and $(P(x)y)^{-1}=P(x^{-1})y^{-1}.$
Furthermore, the fundamental formula 
\begin{equation}\label{eq:fundamental}P(P(x)y)=P(x)P(y)P(x)
\end{equation}
holds for any elements $x,y\in V$, see Proposition II.3.2 (iii) 
~\cite{FK}.

A {\em Euclidean} Jordan algebra is a finite--dimensional real
Jordan algebra $V$ equipped with  an \textit{associative} inner
product $\langle \cdot|\cdot\rangle$, i.e., satisfying $\langle
xy|z\rangle=\langle y|xz\rangle $ for all $x,y,z\in V$. The space
$\Sym$ of $n\times n$ real symmetric
matrices is a Euclidean Jordan algebra with Jordan product
$(1/2)(XY+YX)$ and inner product $\langle
X,Y\rangle={\mathrm{tr}}(XY).$ The spectral theorem (see Theorem
III.1.2 of \cite{FK}) of a Euclidean Jordan algebra $V$ states 
that for $x\in V$ there exist a Jordan frame, i.e., a complete system 
of orthogonal primitive idempotents $c_{1},\dots, c_{r}$, where $r$ 
is the rank of $V$, and real numbers
${\lambda}_{1},\dots,{\lambda}_{r}$, the eigenvalues of $x$, such 
that $x=\sum_{i=1}^{r}{\lambda}_{i}c_{i}.$ Due to the the power 
associative property $x^p \cdot x^q=x^{p+q}$, see \cite{FK}, 
the exponential map $\exp:V\to V$ 
\begin{equation*}
\exp:x\mapsto\sum_{n=0}^{\infty}x^{n}/n!
\end{equation*}
is well defined. Likewise as for the special case $V=\Sym$, the 
Jordan algebra exponential map is 
a bijection between its domain of definition $V$ and its image 
$\Omega:=\exp V.$ In fact, $\Omega$ coincides with the interior 
of the set of square elements of $V$, and this is equal to the set of
invertible squares of $V$. A fundamental theorem of Euclidean
Jordan algebras asserts that (i) $\Omega$ is a symmetric cone, and
(ii) every symmetric cone in a real Euclidean 
space arises in this way. In the case of the Jordan algebra
$V=\Sym$ the Jordan algebra exponential is
simply the matrix exponential map, and hence the corresponding 
symmetric cone is $\Omega=\Spd$, the open convex cone of positive
definite $n\times n$ matrices. Irreducible symmetric cones have 
been completely classified, and the remaining cases consist of 
(i) the cones of positive definite Hermitian and Hermitian quaternion 
$n\times n$ matrices, (ii) the Lorentzian cones, and (iii) a 
27--dimensional exceptional cone. General symmetric cones are Cartesian
products of these. The connected component $\Aut(V)_{\circ}$ of the 
identity $\id_V$ in the Jordan algebra automorphism group 
$\Aut(V)$ is a subgroup of $\Orth(\Omega)$. $\Omega$ is irreducible 
if and only if $V$ is simple, and in this case we have 
$\Aut(V)_{\circ}=\Orth(\Omega)$. For all of these statements, see 
\cite{FK} and the references therein.

In the special case $V=\Sym$ the general formula 
$P(x)y=2x(xy)-x^{2}y$ reduces to $P(X)Y=XYX$, where the 
multiplication in the right hand side of this equation is the usual 
matrix multiplication, not the Jordan multiplication. A key tool 
for generalising $\Sym$ to arbitrary Euclidean Jordan algebras is to 
consistently replace expressions of the form $XYX$ by $P(x)y.$
Throughout this section we will assume that $V$ is a simple
Euclidean Jordan algebra with the associative inner product
$\langle x|y\rangle={\mathrm{tr}}(xy)$ and that $\Omega$ is the 
symmetric cone associated with $V$.

The symmetric cone $\Omega$ carries a $G(\Omega)$\nobreakdash--invariant
Riemannian metric defined by $$\gamma_{x}(u,v)=\langle P(x^{-1})u|
v\rangle, x\in \Omega, u,v\in V$$ for which the Jordan inversion
$x\to x^{-1}$ on $\Omega$ is an involutive isometry fixing $e.$
The curve $t\mapsto P(a^{1/2})(P(a^{-1/2})b)^{t}$ is the unique 
geodesic that joins $a$ to $b$ in $\Omega$, and the Riemannian 
distance $\delta(a,b)$ is given by 
$\delta(a,b)=\bigl(\sum_{i=1}^{n}\log^{2}\lambda_{i}\bigr)^{1/2}$, 
where the $\lambda_{i}$ are the eigenvalues of $P(a^{-1/2})b.$
The \textit{geometric mean} $a\#b$ of two elements $a,b\in\Omega$ 
is defined by 
$$a\#b=P(a^{1/2})(P(a^{-1/2})b)^{1/2}.$$
This is the unique midpoint -- or geodesic middle -- of $a$ and $b$ 
for the 
Riemannian distance $\delta.$ The metric $\delta$ is known as a
Bruhat--Tits metric, i.e., a complete metric satisfying the 
semi--parallelogram law, with midpoint $a\#b$. See \cite{LL} for more details.
If $V=\Sym$ then the geometric mean $A\#B$ of positive definite 
matrices $A$ and $B$ is given by 
$A\#B=A^{1/2}(A^{-1/2}BA^{-1/2})^{1/2}A^{1/2}$. The following basic 
properties of geometric means will be useful for our purpose:

\begin{proposition}[\cite{Lim1}]\label{P:geo}
Let $a$ and $b$ be elements of $\Omega$. Then
\begin{itemize}
\item[(a)] $a\#b$ is  a  unique solution belonging to $\Omega$ of 
the quadratic equation $P(x)a^{-1}=b$.
\item[(b)](The commutativity property) $a\#b=b\#a.$
\item[(c)] (The inversion property)  $(a\#b)^{-1}=a^{-1}\#b^{-1}.$
\item[(d)] $P(a\#b)=P(a)\#P(b)=P(a^{1/2})\left(P(a^{-1/2})P(b)
P(a^{-1/2})\right)^{1/2}P(a^{1/2}).$
\item[(e)] (The transformation property) $g(a\#b)=g(a)\# g(b)$ for 
all $g\in G(\Omega).$
\end{itemize}
\end{proposition}
Let $F(x)=-{\mathrm{ln\ det}}(x)$ be the standard logarithmic
barrier functional on the symmetric cone $\Omega$, where
${\mathrm{det}}$ is the determinant function of the Jordan algebra
$V$, see \cite{FK}. Then one can see that $F'(x)=-x^{-1}$ and the 
Hessian of $F$ is given by $$F''(x)=P(x^{-1}).$$
Proposition \ref{P:geo} implies that the geometric mean
$a\#b^{-1}$ is the scaling point of $a$ and $b\in\Omega$ defined by 
$F$. Indeed,
$F''(a\#b^{-1})a=P\bigl((a\#b^{-1})^{-1}\bigr)a=P(a^{-1}\#b)a=b$,
that is, $w_{F}(a,b)=a\#b^{-1}$.

\begin{lemma}\label{L:key} Let $\Omega$ be an irreducible symmetric 
cone and $\alpha:\Omega\to\Omega$ a function such that 
\begin{equation}\label{eq:condition}
x^{-1}\#y=\alpha(x)^{-1}\#\alpha(y) 
\end{equation}
for all $x,y\in \Omega.$ Then $\alpha=\lambda\cdot\id_{\Omega}$ 
for some positive real number $\lambda.$
\end{lemma}

\begin{proof}
Upon exchanging the roles of $x$ and $y$, Proposition~\ref{P:geo}(a)
implies that Condition~\eqref{eq:condition} 
is equivalent to 
\begin{equation}\label{eq:alphabeta}
\alpha(x)=P(y^{-1}\#x)\alpha(y). 
\end{equation}
Setting $y=e$ and using $e^{-1}\#x=x^{1/2}$ in \eqref{eq:alphabeta}, 
we get 
\begin{equation}\label{eq:firstparagraph}
\alpha(x)=P(x^{1/2})\alpha(e) 
\end{equation}
for all $x\in\Omega$. Let us show that $k(\alpha(e))=\alpha(e)$ for 
all $k\in\Aut(V)_{\circ}$. Applying \eqref{eq:alphabeta} and 
\eqref{eq:firstparagraph} both to $x$ and $y$ we get 
\begin{equation*}
P(x^{1/2})\alpha(e)=P(y^{-1}\#x)\alpha(y)
=P(y^{-1}\#x)\bigl(P(y^{1/2})\alpha(e)\bigr),
\end{equation*}
and hence we obtain the identity 
$\alpha(e)=\bigl(P(x^{-1/2})P(y^{-1}\#x)P(y^{1/2})\bigr)\alpha(e)$ for
all $x,y\in \Omega$, which generalises \eqref{eq:star}. Set 
$$K:=\{P(x^{-1/2})P(y^{-1}\#x)P(y^{1/2}):x,y\in \Omega\}.$$
It follows from the definition of the geometric mean and from the 
fundamental formula that 
\begin{equation*}
P(a\#b)=P(a^{1/2})\Bigl(P\bigl(P(a^{-1/2})b\bigr)\Bigr)^{1/2}P(a^{1/2}).
\end{equation*}
Together with Proposition \ref{P:geo} (b) this implies 
\begin{align*}
P(x^{-1/2})P(y^{-1}\#x)P(y^{1/2})&=P(x^{-1/2})P(x\#y^{-1})
P(y^{1/2})\\
&=P(x^{-1/2})\Bigl(P(x^{1/2})P\bigl(P(x^{-1/2})y^{-1}\bigr)^{1/2}
P(x^{1/2})\Bigr)P(y^{1/2})\\
&=P\bigl(P(x^{1/2})y\bigr)^{-1/2}P(x^{1/2})P(y^{1/2}).
\end{align*}
Therefore, the set $K$ can be written as
$K=\bigl\{P\bigl(P(x)y^{2}\bigr)^{-1/2}P(x)P(y)|x,y\in\Omega\bigr\}$. 
By Koecher's Theorem 4.9 (b)~\cite{Koe}, $K$ generates
$\Aut(V)_{\circ}.$ This implies that the point
$\alpha(e)$ is fixed by all Jordan automorphisms $k\in\Aut(V)_{\circ}.$

Finally, Corollary IV.2.7 \cite{FK} (in which the assumption of 
irreducibility for $\Omega$ is essential) says that the group 
$\Aut(V)_{\circ}$ acts transitively on the set of primitive
idempotents of $V$. The spectral theorem applied to $\alpha(e)$ 
therefore implies that $\alpha(e)=\lambda e$ for some positive real 
number $\lambda$. Together with \eqref{eq:firstparagraph} this 
implies that $\alpha(x)=P(x^{1/2})(\lambda e)=\lambda P(x^{1/2})(e)
=\lambda x$ for all $x\in \Omega$.
\end{proof}

\begin{corollary}\label{C:barrier} Let $H$ be an arbitrary 
self--scaled barrier for the irreducible symmetric cone $\Omega$.
Then there exists a positive constant $\lambda$ such that
$H''(x)=\lambda\cdot F''(x)$ for all $x\in \Omega.$
\end{corollary}
\begin{proof}
Since the Hessian $H''(x)$ are positive definite cone
automorphisms, Lemma \ref{self-adjoint} implies that there exists 
a well--defined function $\Upsilon:\Omega\to\Omega$ such that 
\begin{equation}\label{eq:Ypsilon}
H''(x)=P(\Upsilon(x)^{-1}). 
\end{equation}
Since $H$ is self-scaled, we have 
\begin{align*}
P(\Upsilon(x)^{-1})
&\stackrel{\eqref{eq:Ypsilon},\eqref{eq:thm2ii}}{=}
H''(w_{H})\circ H''(y)\circ H''(w_{H})\\
&=P(\Upsilon(w_{H})^{-1})\circ P(\Upsilon(y)^{-1})
\circ P(\Upsilon(w_{H})^{-1})\\
&\stackrel{\eqref{eq:fundamental}}
=P\Big(P(\Upsilon(w_{H})^{-1})\Upsilon(y)^{-1}\Big)
\end{align*}
for all $x,y\in \Omega,$ where $w_{H}=w_{H}(x,y)$ denotes the 
scaling point of $x$ and $y$ for the self--scaled barrier 
$H$. The quadratic representation $P$ is injective on $\Omega$, 
see Lemma 2.3 \cite{Lim1}. Therefore, the identity above shows that 
\begin{equation*}
\Upsilon(x)^{-1}=P\bigl(\Upsilon(w_{H})^{-1}\bigr)\Upsilon(y)^{-1}
\end{equation*}
for all $x,y\in \Omega.$ By Proposition~\ref{P:geo}, we have 
\begin{equation}\label{eq:Upsilon}
\Upsilon(w_{H})^{-1}=\Upsilon(y)\#\Upsilon(x)^{-1}=\Upsilon(x)^{-1}\#
\Upsilon(y)
\end{equation}
for all $x,y\in \Omega$. Now,
$y=H''(w_{H})(x)=P(\Upsilon(w_{H})^{-1})(x)$ by definition of $w_{H}$,
and Proposition \ref{P:geo} (a) shows that we have
$\Upsilon(w_{H})^{-1}=x^{-1}\#y$, which together with
\eqref{eq:Upsilon} shows that 
$$x^{-1}\#y=\Upsilon(x)^{-1}\#\Upsilon(y)$$ for all 
$x,y\in \Omega.$ The proof is now completed by Lemma \ref{L:key}.
\end{proof}

\begin{theorem}
If $H$ is a self--scaled barrier functional for $\Omega$ then there
exist constants $c_{1}>0$ and $c_{0}\in \Bbb R$ such that
$$H:x\to -c_{1}{\mathrm{ln \ \det}}(x)+c_{0}, \qquad \forall x\in\Omega.$$
\end{theorem}
\begin{proof} Similar to that of Theorem \ref{thm:mainSDP}.
\end{proof}

\section{Acknowledgments}
Both authors would like to thank Professor Leonid Faybusovich 
for suggesting this collaboration and establishing the contact between
them. Raphael Hauser also wishes to express his warmest thanks to 
Syvert N{\o}rsett for inviting him to NTNU Trondheim, and to Peter and
Jonathan Borwein for inviting him to SFU in Vancouver, where part of 
this research was done.

\end{document}